\documentstyle{article}

\def\esp{\raisebox{.3ex}{!`}}

\def\pr{\mbox{\small \O}}

\def\Rel{{\it Rel}}
\def\Fun{{\it Fun}}

\def\LL{${\cal L}^L$}

\def\Mc{${\cal M}_c$}
\def\Lc{${\cal L}_c$}

\def\LLS{${\cal L}^L{\cal S}$}
\def\LRS{${\cal L}^R{\cal S}$}
\def\LS{${\cal LS}$}
\def\LcS{${\cal L}_c{\cal S}$}
\def\CS{${\cal CS}$}
\def\Lcmu{${\cal L}_c\mu$}
\def\DS{${\cal DS}$}
\def\LLSco{${\cal L}^L{\cal S}^{co}$}
\def\MSco{${\cal MS}^{co}$}
\def\McSco{${\cal M}_c{\cal S}^{co}$}
\def\CSco{${\cal CS}^{co}$}
\def\DSco{${\cal DS}^{co}$}

\def\mj{{\mathbf{1}}}
\def\pl{\!+\!}
\def\mn{\!-\!}
\def\cirk{\,{\raisebox{.3ex}{\tiny $\circ$}}\,}
\def\prop#1#2{\vspace{2ex} \noindent{\sc #1.} {\it #2} \par \vspace{2ex}}
\def\dkz{\noindent{\sc Proof. }}
\def\qed{\hfill $\dashv$}
\def\str{\rightarrow}

\begin{document}

\title{Coherence for Monoidal Monads and Comonads}
\author{\small {\sc Kosta Do\v sen} and {\sc Zoran Petri\' c}
\\[1ex]
{\small Mathematical Institute, SANU}\\[-.5ex]
{\small Knez Mihailova 36, p.f.\ 367, 11001 Belgrade,
Serbia}\\[-.5ex]
{\small email: \{kosta, zpetric\}@mi.sanu.ac.rs}}
\date{}
\maketitle

\begin{abstract}
\noindent The goal of this paper is to prove coherence results
with respect to relational graphs for monoidal monads and
comonads, i.e.\ monads and comonads in a monoidal category such
that the endo\-functor of the monad or comonad is a monoidal
functor (this means that it preserves the monoidal structure up to
a natural transformation that need not be an isomorphism). These
results are proved first in the absence of symmetry in the
monoidal structure, and then with this symmetry. The monoidal
structure is also allowed to be given with finite products or
finite coproducts. Monoidal comonads with finite products
axiomatize a plausible notion of identity of deductions in a
fragment of the modal logic S4.
\end{abstract}


\vspace{.5ex}


\section{Introduction}
A monoidal monad is a monad in a monoidal category such that the
endo\-functor of the monad is a monoidal functor, which means that
it preserves the monoidal structure up to a natural transformation
that need not be an isomorphism. (The notion of monoidal functor
stems from \cite{EK66}, Section II.1, and the notion of monoidal
monad from \cite{Ko70} and \cite{Ko72}; for historical remarks on
the notions of monad and comonad see \cite{ML98}, notes at the end
of Chapter VI.) This natural transformation has as components the
arrows
\[
\psi_{A,B}\!:TA\otimes TB\str T(A\otimes B),
\]
and we also have the arrow $\psi_0\!:I\str TI$, which coincides
with the unit arrow $\eta_I$ of the monad; these arrows satisfy
the equations given in Section~3 below.

Our goal in this paper is to prove coherence for this and related
notions of monad with respect to relational graphs; i.e.\ with
respect to the category \Rel, whose arrows are relations between
finite ordinals (and not between any sets, as one would expect by
relying on what is perhaps a more common denotation for \Rel; our
category \Rel\ is the skeleton of the category of relations
between finite sets, but a notation like $Sk(Rel_{fin})$ would be
too cumbersome). Sometimes it will be sufficient to have a
subcategory of \Rel, like the category \Fun\ whose arrows are
functional relations between finite ordinals, or the simplicial
category $\Delta$, whose arrows are order-preserving functional
relations between finite ordinals. Coherence states that there is
a faithful functor from a freely generated monoidal monad, or a
related categorial structure for which we prove coherence, into
\Rel\ or a subcategory of it. We obtain thereby a characterization
of the freely generated monoidal monad, or related categorial
structure, in terms of graphs. Such coherence results give very
useful procedures for deciding whether a diagram of canonical
arrows commutes. (A general treatment of coherence in this spirit
may be found in \cite{DP04}.)

Before we deal with monoidal monads, we consider coherence for
notions of strong monad, which stem from \cite{Ko70}, \cite{Ko72}
and \cite{Mog91}. For monoidal monads we prove coherence first in
the absence of symmetry in the monoidal structure, and then with
this symmetry. After that we prove coherence for monoidal monads
where the monoidal structure is cartesian (i.e.\ with finite
products), or cocartesian (i.e.\ with finite coproducts). This
makes the first part of the paper (Sections 2-6).

A monoidal comonad is a comonad in a monoidal category such that
the endo\-functor of the comonad is a monoidal functor. This
notion is parallel, but not dual, to the notion of monoidal monad.
A dual notion would be the notion of comonoidal comonad, where
instead of $\psi$ and $\psi_0$ we have arrows oriented in the
opposite direction. More recent papers on the notion of monoidal
comonad, or the dual notion of comonoidal monad, called also Hopf
monad, opmonoidal monad or bimonad, are \cite{Moe02}, \cite{MC02},
\cite{BV07} and \cite{PS09}. In \cite{W08}, \cite{BV08} and
references therein, one finds for the canonical arrows of a
comonoidal monad graphical interpretations different from ours and
more involved; coherence is not proved however.

The second part of the paper (Sections 7-11) is parallel to the
first part; instead of monads we have comonads in the same
monoidal contexts, and we prove coherence for these notions. These
results are not dual to those in the first part. By duality, we
can obtain from the results of both parts of the paper coherence
results for comonoidal monads and comonads.

We rely for this paper on \cite{DP09a}, where basic coherence
results for monoidal endo\-functors are proved. Here these
endo\-functors become endo\-functors of monads and comonads. We
presuppose the reader is acquainted with the terminology and
notation of this previous paper, but to make the exposition here
more self-contained we will repeat some definitions.

\section{Coherence for strong monads}
Let a \emph{left monoidal} endo\-functor of a monoidal category
$\langle {\cal A},\otimes,I,a,l,r\rangle$ (in the notation of
\cite{EK66}, Section II.1) be a functor $T$ from $\cal A$ to $\cal
A$ such that the monoidal structure of $\cal A$ is preserved
\emph{locally} by $T$ up to a natural transformation whose
components are the arrows
\[
\psi^L_{A,B}\!:TA\otimes B\str T(A\otimes B);
\]
this means that the following equations hold:
\begin{tabbing}
\hspace{1.5em}\=${(\psi^L a)}$\hspace{3em}\=$Ta_{A,B,C}\cirk
\psi^L_{A\otimes
B,C}\cirk(\psi^L_{A,B}\otimes\mj_C)=\psi^L_{A,B\otimes C}\cirk
a_{TA,B,C}$,
\\[2ex]
\>${(\psi^L r)}$\>$Tr_A\cirk\psi^L_{A,I}=r_{TA}$.
\end{tabbing}

A \emph{right monoidal} endo\-functor is defined analogously with
respect to a natural transformation whose components are the
arrows
\[
\psi^R_{A,B}\!:A\otimes TB\str T(A\otimes B).
\]
The equations corresponding to ${(\psi^L a)}$ and ${(\psi^L r)}$
are
\begin{tabbing}
\hspace{1.5em}\=${(\psi^R a)}$\hspace{3em}\=$Ta_{A,B,C}\cirk
\psi^R_{A\otimes B,C}=\psi^R_{A,B\otimes
C}\cirk(\mj_A\otimes\psi^R_{B,C})\cirk a_{A,B,TC}$,
\\[2ex]
\>${(\psi^R l)}$\>$Tl_A\cirk\psi^R_{I,A}=l_{TA}$.
\end{tabbing}

A \emph{left strong} monad in a monoidal category $\cal A$ is a
monad $\langle T,\eta,\mu\rangle$ (in the notation of \cite{ML98},
Section VI.1) in $\cal A$ such that $T$ is a left monoidal
functor, and we have moreover the equations
\begin{tabbing}
\hspace{1.5em}\=${(\psi^L
\eta)}$\hspace{3em}\=$\psi^L_{A,B}\cirk(\eta_A\otimes\mj_B)\,$\=$=\eta_{A\otimes
B}$,
\\[2ex]
\>${(\psi^L\mu)}$\>$\psi^L_{A,B}\cirk(\mu_A\otimes\mj_B)$\>$=\mu_{A\otimes
B}\cirk T\psi^L_{A,B}\cirk\psi^L_{TA,B}$.
\end{tabbing}
(These equations might be interpreted as saying that $\eta$ is a
\emph{left monoidal} natural transformation from the identity
functor, which is left monoidal, to the left monoidal functor $T$,
while $\mu$ is a left monoidal natural transformation from the
left monoidal functor $TT$ to $T$.)

A \emph{right strong} monad is defined analogously with a right
monoidal functor $T$. The equations corresponding to
${(\psi^L\eta)}$ and ${(\psi^L\mu)}$ are
\begin{tabbing}
\hspace{1.5em}\=${(\psi^R
\eta)}$\hspace{3em}\=$\psi^R_{A,B}\cirk(\mj_A\otimes\eta_B)\,$\=$=\eta_{A\otimes
B}$,
\\[2ex]
\>${(\psi^R
\mu)}$\>$\psi^R_{A,B}\cirk(\mj_A\otimes\mu_B)$\>$=\mu_{A\otimes
B}\cirk T\psi^R_{A,B}\cirk\psi^R_{A,TB}$.
\end{tabbing}
The notions of left strong and right strong monad are derived from
\cite{Ko70} and \cite{Ko72} in \cite{Mog91}.

Let \LLS\ be the category of the left strong monad freely
generated by an arbitrary set of objects, and let $\Delta$ be the
simplicial category (see \cite{ML98}, Section VII.5); the arrows
of $\Delta$ are order-preserving functions between finite
ordinals. We define a functor $G$ from \LLS\ to $\Delta$ by
stipulating that $GA$ for $A$ an object \LLS\ is the number of
occurrences of $T$ in $A$, while $Gf$ for an arrow term $f$ of
\LLS\ is defined inductively on the complexity of $f$. If $f$ is
$a_{A,B,C}$, $l_A$, $r_A$, or $\psi^L_{A,B}$, then $Gf$ is the
identity function; next, we have clauses corresponding to the
following pictures:
\begin{center}
\begin{picture}(280,40)

\put(70,10){\circle*{2}} \put(102,10){\circle*{2}}
\put(54,10){\circle*{2}} \put(70,30){\circle*{2}}
\put(102,30){\circle*{2}}

\put(80,20){$\ldots$}

\put(70,12){\line(0,1){16}} \put(102,12){\line(0,1){16}}

\put(70,7){\makebox(0,0)[t]{\tiny$G\!A\!\mn\! 1$}}
\put(102.5,7){\makebox(0,0)[t]{\tiny$0$}}
\put(54.5,7){\makebox(0,0)[t]{\tiny$G\!A$}}
\put(70,33){\makebox(0,0)[b]{\tiny$G\!A\!\mn\! 1$}}
\put(102.5,33){\makebox(0,0)[b]{\tiny$0$}}
\put(35,20){\makebox(0,0)[r]{$G\eta_A$}}

\put(210,10){\circle*{2}} \put(257,10){\circle*{2}}
\put(194,30){\circle*{2}} \put(210,30){\circle*{2}}
\put(257,30){\circle*{2}}

\put(235,20){$\ldots$}

\put(225.5,7){\makebox(0,0)[t]{\tiny$G\!A\!\mn\! 1$}}

\put(225.5,33){\makebox(0,0)[b]{\tiny$G\!A\!\mn\! 1$}}

\put(225,12){\line(0,1){16}}

\put(225,10){\circle*{2}}

\put(225,30){\circle*{2}}

\put(210,12){\line(0,1){16}} \put(257,12){\line(0,1){16}}
\put(210,12){\line(-1,1){16}}

\put(210,7){\makebox(0,0)[t]{\tiny$G\!A$}}
\put(257.5,7){\makebox(0,0)[t]{\tiny$0$}}
\put(194,33){\makebox(0,0)[b]{\tiny$G\!A\!\pl\! 1$}}
\put(210,33){\makebox(0,0)[b]{\tiny$G\!A$}}
\put(257.5,33){\makebox(0,0)[b]{\tiny$0$}}
\put(175,20){\makebox(0,0)[r]{$G\mu_A$}}

\end{picture}
\end{center}
\begin{center}
\begin{picture}(160,30)(0,10)

\put(10,16){\line(0,1){20}} \put(30,16){\line(0,1){20}}

\put(10,16){\line(1,0){20}} \put(10,36){\line(1,0){20}}

\put(40,16){\line(0,1){20}} \put(60,16){\line(0,1){20}}

\put(40,16){\line(1,0){20}} \put(40,36){\line(1,0){20}}

\put(14,25){\makebox(0,0)[l]{$Gf$}}
\put(44,25){\makebox(0,0)[l]{$Gg$}}

\put(-45,25){\makebox(0,0)[l]{$G(f\otimes g)$}}

\put(150,16){\line(0,1){20}} \put(170,16){\line(0,1){20}}

\put(150,16){\line(1,0){20}} \put(150,36){\line(1,0){20}}

\put(140,16){\line(0,1){20}}

\put(154,25){\makebox(0,0)[l]{$Gf$}}

\put(95,25){\makebox(0,0)[l]{$GT\!f$}}

\end{picture}
\end{center}
and  $G(g\cirk f)$ is the composition of the functions $Gf$ and
$Gg$. We can prove the following.

\prop{\LLS-Coherence}{The functor $G$ from \LLS\ to $\Delta$ is
faithful.}

\dkz By naturality and functorial equations, and the equations
${(\psi^L\eta)}$ and ${(\psi^L\mu)}$, every arrow term of
$f\!:A\str B$ of \LLS\ is equal to an arrow term $f_2\cirk
f_1\!:A\str B$ such that $\eta$ and $\mu$ do not occur in
$f_1\!:A\str C$, and $a$, $l$, $r$ and $\psi^L$ do not occur in
$f_2\!:C\str B$. We can uniquely determine $Gf_2$ from $Gf$ (since
$Gf_1$ is the identity function, $Gf_2=Gf$), and from $B$ and
$Gf_2$ we obtain a unique $C$ as a possible source of $f_2$. By
the isomorphism of $\Delta$ with the monad freely generated by a
single object (see the references in \cite{DP08}, Section~3), we
can conclude that $f_2$ is uniquely determined by $Gf$, while
$f_1$ is unique by \LL-Coherence (see \cite{DP09a}, Section~5).
\qed

\vspace{2ex}

Let \LRS\ be the category of the right strong monad freely
generated by an arbitrary set of objects. It is easy to show that
this category is isomorphic to \LLS; it is a mirror image of \LLS.
Hence \LRS-Coherence could be a result, exactly analogous to
\LLS-Coherence, about a faithful functor from \LRS\ to the
simplicial category $\Delta$. For future use, however, we need
another coherence result for \LRS, which is given with respect to
a functor $G$ from \LRS\ to the category \Fun\ whose arrows are
arbitrary functions between finite ordinals. This functor is
defined as $G$ from \LLS\ to $\Delta$ except for the additional
clause for $G\psi^R_{A,B}$ corresponding to the following picture:
\begin{center}
\begin{picture}(320,50)

\put(142,16){\line(0,1){20}} \put(146,16){\line(0,1){20}}

\put(141.75,16){\line(1,0){4.3}} \put(141.75,36){\line(1,0){4.3}}

\put(170,16){\line(0,1){20}} \put(174,16){\line(0,1){20}}

\put(169.75,16){\line(1,0){4.3}} \put(169.75,36){\line(1,0){4.3}}

\put(132,16){\line(3,2){30}}

\put(130,42){\makebox(0,0)[l]{$\;\;\;\;A\otimes T\,B$}}
\put(128.25,8){\makebox(0,0)[l]{$T\,(A\;\otimes\;\,B)$}}

\end{picture}
\end{center}
In this picture we obtain a crossing if there is a $T$ in $A$.
Hence with this clause we abandon the category $\Delta$, and must
consider also functions that are not order-preserving. For this
functor G we can prove the following.

\prop{\LRS-Coherence}{The functor $G$ from \LRS\ to \Fun\ is
faithful.}

The proof is as for \LLS-Coherence, except for the parenthetical
remark about determining $Gf_2$ out of $Gf$. Now $Gf_2$ is not
equal to $Gf$, but it is still uniquely determined by it. Since
$Gf_1$ is a bijection, and $Gf_2$ is an order-preserving function,
the latter is uniquely determined by knowing for each element of
the target of $Gf$ the number of elements of the source of $Gf$
mapped to this element by $Gf$. Note that we could prove
\LLS-Coherence with respect to a functor from \LLS\ to \Fun\
analogous to the functor $G$ from \LRS\ to \Fun.

\section{Coherence for monoidal monads}
Let a \emph{locally} monoidal endo\-functor of a monoidal category
$\langle {\cal A},\otimes,I,a,l,r\rangle$ be a functor $T$ that is
both left monoidal and right monoidal, and we have moreover the
equation
\begin{tabbing}
${(\psi^L\psi^R a)}$\hspace{1.4em}$Ta_{A,B,C}\cirk\psi^L_{A\otimes
B,C}\cirk(\psi^R_{A,B}\otimes\mj_C)= \psi^R_{A,B\otimes
C}\cirk(\mj_A\otimes\psi^L_{B,C})\cirk a_{A,TB,C}$.
\end{tabbing}

A \emph{monoidal monad} in a monoidal category $\cal A$ is a monad
$\langle T,\eta,\mu\rangle$ in $\cal A$ both left strong and right
strong, such that $T$ is a locally monoidal endo\-functor, and we
have moreover the equation
\begin{tabbing}
\hspace{1.5em}\=${(\psi a)}$\hspace{3em}\=$Ta_{A,B,C}\cirk
\psi_{A\otimes B,C}\cirk(\psi_{A,B}\otimes\mj_{TC})=$\kill

\>${(\psi^L\psi^R\mu)}$\>$\mu_{A\otimes B}\cirk
T\psi^L_{A,B}\cirk\psi^R_{TA,B}=\mu_{A\otimes B}\cirk
T\psi^R_{A,B}\cirk\psi^L_{A,TB}$:\\*[1ex] \`$TA\otimes TB\str
T(A\otimes B)$.
\end{tabbing}

An alternative definition of monoidal monad is obtained by
stipulating that in a monoidal category $\cal A$ we have a monad
$\langle T,\eta,\mu\rangle$ and a natural transformation whose
components are the arrows
\[
\psi_{A,B}\!:TA\otimes TB\str T(A\otimes B),
\]
which satisfy the equations
\begin{tabbing}
\hspace{1.5em}\=${(\psi a)}$\hspace{3em}\=$Ta_{A,B,C}\cirk
\psi_{A\otimes B,C}\cirk(\psi_{A,B}\otimes\mj_{TC})=$\kill

${(\psi a)}$\hspace{1.2em}\=$Ta_{A,B,C}\cirk \psi_{A\otimes
B,C}\cirk(\psi_{A,B}\otimes\mj_{TC})=\psi_{A,B\otimes
C}\cirk(\mj_{TA}\otimes\psi_{B,C})\cirk a_{TA,TB,TC}$,
\\[2ex]
\>${(\psi
l)}$\>$Tl_A\cirk\psi_{I,A}\cirk(\eta_I\otimes\mj_{TA})\;$\=$=l_{TA}$,
\\[2ex]
\>${(\psi
r)}$\>$Tr_A\cirk\psi_{A,I}\cirk(\mj_{TA}\otimes\eta_I)$\>$=r_{TA}$,
\\[3ex]
\>${(\psi
\eta)}$\>$\psi_{A,B}\cirk(\eta_A\otimes\eta_B)\;$\=$=\eta_{A\otimes
B}$,\\*[2ex] \>${(\psi
\mu)}$\>$\psi_{A,B}\cirk(\mu_A\otimes\mu_B)$\>$=\mu_{A\otimes
B}\cirk T\psi_{A,B}\cirk\psi_{TA,TB}$.
\end{tabbing}
The first three of these equations, together with $\psi_0=\eta_I$,
say that $T$ is a monoidal functor, while the last two equations,
together with $\psi_0=\eta_I$ and the monad equation $\mu_I\cirk
T\eta_I=\mj_{TI}$, ensure that $\eta$ and $\mu$ are monoidal
natural transformations in the sense of \cite{EK66} (Section II.1;
see also \cite{ML98}, Section XI.2). Coherence for monoidal
endo\-functors is proved in \cite{DP09a} (Section~4).

With $\psi_{A,B}$ being defined as either of the two sides of the
equation ${(\psi^L\psi^R\mu)}$, and with
\begin{tabbing}
\hspace{1.5em}\=${(\psi a)}$\hspace{3em}\=$Ta_{A,B,C}\cirk
\psi_{A\otimes B,C}\cirk(\psi_{A,B}\otimes\mj_{TC})=$\kill

\>\>$\psi^L_{A,B}\,$\=$=_{df}\psi_{A,B}\cirk(\mj_{TA}\otimes\eta_B)$,\\*[2ex]
\>\>$\psi^R_{A,B}$\>$=_{df}\psi_{A,B}\cirk(\eta_A\otimes\mj_{TB})$,
\end{tabbing}
we can show that the two definitions of monoidal monad amount to
the same notion. Both of these definitions stem from \cite{Ko70}
and \cite{Ko72}.

Let \LS\ be the category of the monoidal monad, with $\psi^L$ and
$\psi^R$ primitive, freely generated by an arbitrary set of
objects. We define the functor $G$ from \LS\ to the category \Fun\
by combining what we had in the preceding section for the functors
$G$ from \LLS\ to $\Delta$ and \LRS\ to \Fun. We can prove the
following.

\prop{\LS-Coherence}{The functor $G$ from \LS\ to \Fun\ is
faithful.}

\dkz By naturality and functorial equations, and the equations
${(\psi^L\eta)}$, ${(\psi^L\mu)}$, ${(\psi^R\eta)}$ and
${(\psi^R\mu)}$, every arrow term of $f\!:A\str B$ of \LS\ is
equal to an arrow term $f_2\cirk f_1\!:A\str B$ such that $\eta$
and $\mu$ do not occur in $f_1\!:A\str C$, and $a$, $l$, $r$,
$\psi^L$ and $\psi^R$ do not occur in $f_2\!:C\str B$. We can
uniquely determine $Gf_2$ from $Gf$, and as in the proof of
\LLS-Coherence and \LRS-Coherence, from $B$ and $Gf_2$ we obtain a
unique $C$ as a possible source of $f_2$. As in the previous
proofs, we can conclude that $f_2$ is uniquely determined by $Gf$.
However, $Gf_1$ is not thereby uniquely determined, as the two
sides of ${(\psi^L\psi^R\mu)}$ show. Hence $f_1$ is not unique
either.

Let the normal form of the proof of $\cal L$-Coherence of
\cite{DP09a} (Section~5) for $f_1$ be $g_m\ldots g_1$. If in this
normal form for a $\psi^R$-factor $g_i$ and a $\psi^L$-factor
$g_{i+1}$ we have $\tau(g_{i+1})=\tau(g_i)\pl 1$ and the function
$Gf_2$ has the same value when it is applied to $\tau(g_i)$ and
$\tau(g_{i+1})$, then we rely essentially on the equation
${(\psi^L\psi^R\mu)}$ to permute $g_i$ with $g_{i+1}$. By
proceeding in this manner, we obtain an arrow term $f'_1$ such
that the permutation $Gf'_1$ has the least possible number of
inversions. By $\cal L$-Coherence, $f'_1$ is unique. \qed

\section{Coherence for symmetric monoidal monads}
A \emph{locally linear} endo\-functor is a locally monoidal
endo\-functor $T$ in a symmetric monoidal category that satisfies
the equation
\begin{tabbing}
\hspace{1.5em}\=${(\psi c)}$\hspace{3em}\=$Tc_{A,B}\cirk
\psi_{A,B}=\psi_{B,A}\cirk c_{TA,TB}$.\kill

\>${(\psi^L\psi^R
c)}$\>$Tc_{A,B}\cirk\psi^L_{A,B}=\psi^R_{B,A}\cirk c_{TA,B}$.
\end{tabbing}
From this equation, which says that $T$ preserves $c$ locally up
to $\psi^L$ and $\psi^R$, we obtain immediately a definition of
$\psi^R$ in terms of $\psi^L$, and vice versa.

A \emph{symmetric} monoidal monad is a monoidal monad in a
symmetric monoidal category whose endo\-functor is a locally
linear monoidal endo\-functor. In the language of monoidal monads
with $\psi$ primitive, the equation ${(\psi^L\psi^R c)}$ of
symmetric monoidal monads is replaced by the equation
\begin{tabbing}
\hspace{1.5em}\=${(\psi c)}$\hspace{3em}\=$Tc_{A,B}\cirk
\psi_{A,B}=\psi_{B,A}\cirk c_{TA,TB}$.\kill

\>${(\psi c)}$\>$Tc_{A,B}\cirk\psi_{A,B}=\psi_{B,A}\cirk
c_{TA,TB}$
\end{tabbing}
(see \cite{DP09a}, Section~6).

Let \LcS\ be the category of the symmetric monoidal monad, with
$\psi^L$ and $\psi^R$ primitive, freely generated by an arbitrary
set of objects. We define the functor $G$ from \LcS\ to the
category \Fun\ by stipulating first that $GA$ is the number of
occurrences of generating objects, i.e.\ propositional letters, in
$A$ plus the number of occurrences of $T$ in $A$. Up to now we
took $GA$ to be just the number of occurrences of $T$ in $A$, but
we could as well have counted also occurrences of generating
objects; this was however superfluous up to now. The remainder of
the definition of $G$ is as the definition of the functor $G$ from
\LS\ to \Fun, with the additional standard clause for $Gc_{A,B}$
corresponding to the picture
\begin{center}
\begin{picture}(30,50)

\put(2,16){\line(4,5){16}}

\put(6,16){\line(4,5){16}}

\put(2,16){\line(1,0){4}}

\put(18,36){\line(1,0){4}}

\put(18,16){\line(-4,5){16}}

\put(22,16){\line(-4,5){16}}

\put(2,36){\line(1,0){4}}

\put(18,16){\line(1,0){4}}

\put(0,42){\makebox(0,0)[l]{$A\otimes B$}}
\put(-1.75,8){\makebox(0,0)[l]{$B\otimes A$}}

\end{picture}
\end{center}
We can prove the following.

\prop{\LcS-Coherence}{The functor $G$ from \LcS\ to \Fun\ is
faithful.}

\noindent For the proof we proceed as for \LS-Coherence, by
relying on \Lc-Coherence of \cite{DP09a} (Section~6).

\section{Coherence for cartesian monoidal monads}
A \emph{cartesian} monoidal monad (not to be confused with the
cartesian monads of \cite{Lei03}, Section 4.1) is a symmetric
monoidal monad in a cartesian category (by which is meant a
monoidal category whose monoidal structure is given by finite
products), which satisfies moreover the equation
\begin{tabbing}
\hspace{1.5em}\=${(\psi c)}$\hspace{3em}\=$Tc_{A,B}\cirk
\psi_{A,B}=\psi_{B,A}\cirk c_{TA,TB}$.\kill

\>${(\psi\Delta)}$\>$T\Delta_A=\psi_{A,A}\cirk\Delta_{TA}$,
\end{tabbing}
where $\Delta_A\!:A\str A\otimes A$ is a component of the diagonal
natural transformation of the cartesian structure. The equation
${(\psi\Delta)}$ says intuitively that $T$ preserves $\Delta$ up
to $\psi$. In the terminology of \cite{DP09a} (Section~7), $T$ is
a conjunctive relevant endo\-functor.

Note that in the definition of cartesian monoidal monad we do not
assume the equation
\begin{tabbing}
\hspace{1.5em}\=${(\psi c)}$\hspace{3em}\=$Tc_{A,B}\cirk
\psi_{A,B}=\psi_{B,A}\cirk c_{TA,TB}$.\kill

\>${(\psi\esp)}$\>$T\esp_A=\eta_I\cirk\esp_{TA}$,
\end{tabbing}
where $\esp_A\!:A\str I$ is the unique arrow from $A$ to the
terminal object (empty product) $I$. The equation ${(\psi\esp)}$
says intuitively that $T$ preserves $\esp$ up to $\psi$, but we
cannot assume that equation if we are guided by coherence, as the
following pictures show:
\begin{center}
\begin{picture}(130,50)

\put(3,15.5){\line(0,1){16}}
\put(0,38){\makebox(0,0)[l]{$TA$}}
\put(0,8){\makebox(0,0)[l]{$TI$}}

\put(-30,23){\makebox(0,0)[l]{$T\esp_A$}}

\put(110,48){\makebox(0,0)[l]{$TA$}}
\put(110,23){\makebox(0,0)[l]{$\;\;\;I$}}
\put(110,-2){\makebox(0,0)[l]{$TI$}}

\put(80,35){\makebox(0,0)[l]{$\esp_{TA}$}}
\put(80,11){\makebox(0,0)[l]{$\eta_I$}}
\end{picture}
\end{center}

We will now deal with coherence for the notion of cartesian
monoidal monad. Let \CS\ be the category of the cartesian monoidal
monad freely generated by an arbitrary set of objects. We define
the functor $G$ from \CS\ to the category \Rel, whose arrows are
arbitrary relations between finite ordinals, by adding to the
definition of $G$ from \LcS\ to \Fun\ the clause for $G\Delta_A$
that corresponds to the following picture:
\begin{center}
\begin{picture}(25,35)(0,7)

\put(4,16){\line(1,3){6.5}}

\put(8,16){\line(1,3){6.5}}

\put(4,16){\line(1,0){4}}

\put(10.5,35.5){\line(1,0){4}}

\put(17,16){\line(-1,3){6.5}}

\put(21,16){\line(-1,3){6.5}}

\put(17,16){\line(1,0){4}}

\put(9,42){\makebox(0,0)[l]{$A$}}
\put(-1.75,8){\makebox(0,0)[l]{$A\otimes A$}}

\end{picture}
\end{center}
and the clause that says that $G\esp_A$ is the empty relation
between $GA$ and \pr\ (see \cite{DP09a}, Sections 7-8).

As an auxiliary result for the proof of \CS-Coherence we establish
first a lemma for which we need the notions of diversified object
and scope of \cite{DP09a} (Section~2). An object that is a
propositional formula is diversified when every generator (which
means generating object, i.e.\ propositional letter, or generating
functor) occurs in it at most once; for $E^iC$ a subformula of
$D$, the scope in $D$ of the outermost occurrence of $E^i$ in
$E^iC$ is the set of all the generators in $C$. The category \Lc\
is the free symmetric monoidal category with a family of locally
linear endo\-functors (see the beginning of the preceding section
and \cite{DP09a}, Section~6). Here is our auxiliary lemma.

\prop{\Lc-Theoremhood Lemma}{For $A$ and $B$ diversified objects,
there is an arrow $f\!:A\str B$ of \Lc\ iff the generators of $A$
and $B$ coincide, and for every generating functor $E^i$ of $A$
the scope of $E^i$ in $A$ is a subset of the scope of $E^i$ in
$B$.}

\dkz We proceed as in the proof of the \Mc-Theoremhood Lemma of
\cite{DP09a} (Section~7), except for the assumption about the form
of $A$ in the induction step, which is now
\[
D_1\otimes A_1\otimes\ldots\otimes D_{i-1}\otimes A_{i-1}\otimes
D_i\otimes E^1A_i\otimes D_{i+1}\otimes
A_{i+1}\otimes\ldots\otimes D_n\otimes A_n\otimes D_{n+1},
\]
with the union of the generators of $A_1,\ldots,A_n$ making the
scope of $E^1$ in $B$.\qed

\vspace{2ex}

Let \Lcmu\ be the category defined like \Lc, with $\psi^l$ and
$\psi^R$ primitive, save that we have in addition the primitive
arrow terms $\mu^i_A\!:E^iE^iA\str E^iA$ for which we assume the
naturality and associativity equations of monads, and the
equations ${(\psi^L\mu)}$ and ${(\psi^L\psi^R\mu)}$ with $\mu$ and
$T$ replaced respectively by $\mu^i$ and $E^i$. This category does
not differ essentially from \LcS\ in which $\eta$ is absent. From
the proof of \LcS-Coherence we may easily infer
\Lcmu-\emph{Coherence} in the following form:
\begin{quote}
{\it For all arrow terms $f,g\!:A\str B$ of \Lcmu\ with $B$
diversified we have $f=g$ in \Lcmu.}
\end{quote}
We define $\psi^i$ in \Lcmu\ as we did previously in Section~3.

We have introduced the category \Lcmu\ to formulate the following
lemma without complications involving graphs.

\prop{\Lcmu-Theoremhood Lemma}{For $A$ diversified on generating
objects and $B$ diversified, there is an arrow $f\!:A\str B$ of
\Lc\ iff the generators of $A$ and $B$ coincide, and for every
generating functor $E^i$ of $B$ the union of the scopes of the
occurrences of $E^i$ in $A$ is a subset of the union of the scope
of $E^i$ in $B$ with the set $\{E^i\}$.}

\dkz The union with the set $\{E^i\}$ is mentioned above because
an occurrence of $E^i$ may be in the scope of another occurrence
of $E^i$ in $A$. Such an occurrence of $E^i$ is called
\emph{nested}.

In the beginning we proceed for this proof as for the proof of the
\Mc-Theoremhood Lemma of \cite{DP09a} (Section~7) until the
assumption about the form of $A$ in the induction step of the main
induction. This form is now like the form mentioned in the proof
of the \Lc-Theoremhood Lemma save that some of the $A_j$'s for
$j\neq i$ may be replaced by $E^1A_j$.

Then we have an auxiliary induction on the number $m$ of nested
occurrences of $E^1$ in $A$, in order to prove that there is an
arrow $f'$ of \Lcmu\ from $A$ to the formula $A'$ obtained from
$A$ by deleting all the nested occurrences of $E^1$. In the basis
of this auxiliary induction, when $m=0$, we have an identity
arrow. In the induction step of this auxiliary induction, since
$E^1$ is not in the scope of any $E$ in $B$, we have a subformula
of $A$ like $E^1(C_1\otimes E^1C_2\otimes C_3)$. Let $A''$ be
obtained from $A$ by replacing this subformula with
$E^1(C_1\otimes C_2\otimes C_3)$. It is clear that we have an
arrow $g\!:A\str A''$ of \Lcmu, and by the induction hypothesis of
the auxiliary induction we have a desired arrow $f''\!:A''\str
A'$. So we have a desired arrow $f'\!:A\str A'$.

Then as in the \Mc-Theoremhood and \Lc-Theoremhood Lemmata we have
a desired arrow from $A'$ to \[E^1(A_1\otimes\ldots\otimes
A_n)\otimes D_1\otimes\ldots\otimes D_{n+1},\] and we proceed as
before for the remainder of the proof. \qed

\vspace{2ex}

We can then prove the following lemma analogous to Lemma~1 of
\cite{DP09a} (Section~7).

\prop{Lemma~1}{For the arrow term \[ f\!:A[EA_1\otimes EA_2]\str
B\] of \Lcmu\ and $g\!:B\str C$ a $\mu$-factor such that the
ordinals corresponding to the outermost occurrences of $E$ in
$EA_1$ and $EA_2$ are respectively $i$ and $j$, and
$(GHf)(i)\neq(GHf)(j)$, while $(GH(g\cirk f))(i)=(GH(g\cirk
f))(j)$, there exists an arrow term \[ f'\!:A[E(A_1\otimes
A_2)]\str C\] of \Lcmu\ such that $g\cirk f=f'\cirk
A[\psi_{A_1,A_2}]$.}

\dkz Note first that every arrow term of \Lcmu\ is a substitution
instance of an arrow term of \Lcmu\ with a diversified target. So
we may assume that $C$ in the lemma is diversified. That $f'$
exists follows from the assumption  that we have $g\cirk f$ and
from the \Lcmu-Theoremhood Lemma. That $g\cirk f=f'\cirk
A[\psi_{A_1,A_2}]$ follows from \Lcmu-Coherence. \qed

\vspace{2ex}

We can now prove the following.

\prop{\CS-Coherence}{The functor $G$ from \CS\ to \Rel\ is
faithful.}

\dkz We establish first that every arrow term $f$ of \CS\ is equal
to an arrow term $f_3\cirk f_2\cirk f_1$ such that in the
developed arrow term $f_1$ all the heads of factors are of the
form $\Delta_A$ or $\esp_A$ with $A$ atomic (see \cite{DP09a},
Section~2, for the notions of developed, head, factor and atomic),
while $f_2$ is an arrow term of \LcS\ without occurrences of
$\eta$, and in the developed arrow term $f_3$ all the heads are of
the form $\eta_B$. This is established as for Lemma~2 of
\cite{DP09a} (Section~7) and for $\cal C$-Coherence of
\cite{DP09a} (Section~8); we rely moreover on naturality and
functorial equations to produce $f_3$.

The remainder of the proof is analogous to the proof of Lemma~5 of
\cite{DP09a} (Section~7), which is based on Lemma~3 (ibid.). The
notion of short circuit is the same, but in the proof of Lemma~3
we have that $g$ is a $\mu$-factor, and not a $\psi$-factor.

That $Gf$ determines uniquely $Gf_1$, $Gf_2$, $Gf_3$, and the
targets of $f_1$ and $f_2$ is established as at the end of the
proof of Proposition~5 of \cite{DP09a} (Section~8). We may assume
that the target of $f$ is $\otimes$-free, and so $f_1$ will be
$\Delta$-free. \qed

\section{Coherence for cocartesian monoidal monads}
A \emph{cocartesian} monoidal monad is a symmetric monoidal monad
in a cocartesian category, by which we mean a monoidal category
whose monoidal structure is given by finite coproducts; we have
moreover the equation
\begin{tabbing}
\hspace{1.5em}\=${(\psi c)}$\hspace{3em}\=$Tc_{A,B}\cirk
\psi_{A,B}=\psi_{B,A}\cirk c_{TA,TB}$.\kill

\>${(\psi\;\mbox{\it def}\,)}$\>$\psi_{A,B}=\nabla_{T(A\otimes
B)}\cirk (T\iota^1_{A,B}\otimes T\iota^2_{A,B})$,
\end{tabbing}
where $\nabla_A\!:A\otimes A\str A$ is a component of the
codiagonal natural transformation, while $\iota^1\!:A\str A\otimes
B$ and $\iota^2\!:B\str A\otimes B$ are components of the
injection natural transformations of the cocartesian structure.

More simply, we can define a cocartesian monoidal monad as a
cocartesian category with a monad in it. The definition of $\psi$
is given by ${(\psi\;\mbox{\it def}\,)}$.

In every cocartesian monoidal monad, the functor $T$ preserves
$\nabla$ up to $\psi$, in the sense that we have the equation
\[
T\nabla_A\cirk\psi_{A,A}=\nabla_{TA}.
\]
If $!_A\!:I\str A$ is the unique arrow from the initial object $I$
to $A$, then $T$ preserves also $!$ up to $\psi_0$, which is
defined as $\eta_I$, in the sense that we have the equation
\[
T!_A\cirk\eta_I=!_{TA}.
\]

We will now deal with coherence for the notion of cocartesian
monoidal monad. Let \DS\ be the category of the cocartesian
monoidal monad freely generated by an arbitrary set of objects. We
define the functor $G$ from \DS\ to the category \Fun, whose
arrows are arbitrary relations between finite ordinals, by adding
to the definition of $G$ from \LcS\ to \Fun\ the clause for
$G\nabla_A$ that corresponds to the following picture:
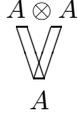
\begin{center}
\begin{picture}(30,35)(0,7)
\put(4,35.5){\line(1,-3){6.5}}

\put(8,35.5){\line(1,-3){6.5}}

\put(4,35.5){\line(1,0){4}}

\put(10.5,16){\line(1,0){4}}

\put(17,35.5){\line(-1,-3){6.5}}

\put(21,35.5){\line(-1,-3){6.5}}

\put(17,35.5){\line(1,0){4}}

\put(9,8){\makebox(0,0)[l]{$A$}}
\put(0,42){\makebox(0,0)[l]{$A\otimes A$}}

\end{picture}
\end{center}
and the clause that says that $G!_A$ is the empty function from
$\pr$ to $GA$ (see \cite{DP09a}, Section~9). Then we can prove the
following.

\prop{\DS-Coherence}{The functor $G$ from \DS\ to \Fun\ is
faithful.}

\dkz We establish first that every arrow term $f$ of \DS\ is equal
to an arrow term $f_2\cirk f_1$ such that $f_1$ is an arrow term
of \LcS, and in the developed arrow term $f_2$ every factor is
either a $\nabla$-factor or a $!$-factor such that the index of
the head is a propositional letter.

Then we can ascertain that the target of $f_1$, which is the
source of $f_2$, is uniquely determined by $Gf$. The graph $Gf_2$
is uniquely determined by $Gf$, and so is the graph $Gf_1$ if in
$f_1$ we get rid of useless crossings (cf.\ \cite{DP09a},
Section~7). Then we rely on \LcS-Coherence of Section~4 and $\cal
D$-Coherence of \cite{DP09a} (Section~9) to obtain \DS-Coherence.
(As a matter of fact, what we need is a rather trivial instance of
$\cal D$-Coherence.) \qed

\section{Coherence for strong comonads}
A \emph{left strong} comonad in a monoidal category $\cal A$ is a
comonad $\langle L,\varepsilon,\delta\rangle$ (in the notation of
\cite{ML98}, Section VI.1) in $\cal A$ such that $L$ is a left
monoidal functor (see Section~2), and we have moreover the
equations
\begin{tabbing}
\hspace{1.5em}\=${(\psi^L
\varepsilon)}$\hspace{3em}\=$\varepsilon_{A\otimes
B}\cirk\psi^L_{A,B}\,$\=$=\varepsilon_A\otimes\mj_B$,
\\[2ex]
\>${(\psi^L\delta)}$\>$\delta_{A\otimes B}\cirk\psi^L_{A,B}$\>$=
L\psi^L_{A,B}\cirk\psi^L_{LA,B}\cirk(\delta_A\otimes\mj_B)$.
\end{tabbing}
(These equations might be interpreted as saying that $\varepsilon$
is a left monoidal natural transformation from $L$ to the identity
functor, while $\delta$ is a left monoidal natural transformation
from $L$ to $LL$.)

A \emph{right strong} comonad is defined analogously with a right
monoidal functor $L$. The equations corresponding to
${(\psi^L\varepsilon)}$ and ${(\psi^L\delta)}$ are
\begin{tabbing}
\hspace{1.5em}\=${(\psi^R
\varepsilon)}$\hspace{3em}\=$\varepsilon_{A\otimes
B}\cirk\psi^R_{A,B}\,$\=$=\mj_A\otimes\varepsilon_B$,
\\[2ex]
\>${(\psi^R\delta)}$\>$\delta_{A\otimes B}\cirk\psi^R_{A,B}$\>$=
L\psi^R_{A,B}\cirk\psi^R_{A,LB}\cirk(\mj_A\otimes\delta_B)$.
\end{tabbing}

Let \LLSco\ be the category of the left strong comonad freely
generated by an arbitrary set of objects. We define a functor $G$
from \LLSco\ to the category $\Delta^{op}$ as the functor $G$ from
\LLS\ to the simplicial category $\Delta$ in Section~2, with the
clauses for $G\eta_A$ and $G\mu_A$ replaced by dual clauses, which
correspond to the following pictures:
\begin{center}
\begin{picture}(280,40)

\put(70,10){\circle*{2}} \put(102,10){\circle*{2}}
\put(54,30){\circle*{2}} \put(70,30){\circle*{2}}
\put(102,30){\circle*{2}}

\put(80,20){$\ldots$}

\put(70,12){\line(0,1){16}} \put(102,12){\line(0,1){16}}

\put(70,7){\makebox(0,0)[t]{\tiny$G\!A\!\mn\! 1$}}
\put(102.5,7){\makebox(0,0)[t]{\tiny$0$}}
\put(54.5,33){\makebox(0,0)[b]{\tiny$G\!A$}}
\put(70,33){\makebox(0,0)[b]{\tiny$G\!A\!\mn\! 1$}}
\put(102.5,33){\makebox(0,0)[b]{\tiny$0$}}
\put(35,20){\makebox(0,0)[r]{$G\varepsilon_A$}}

\put(210,10){\circle*{2}} \put(257,10){\circle*{2}}
\put(194,10){\circle*{2}} \put(210,30){\circle*{2}}
\put(257,30){\circle*{2}}

\put(235,20){$\ldots$}

\put(225.5,7){\makebox(0,0)[t]{\tiny$G\!A\!\mn\! 1$}}

\put(225.5,33){\makebox(0,0)[b]{\tiny$G\!A\!\mn\! 1$}}

\put(225,12){\line(0,1){16}}

\put(225,10){\circle*{2}}

\put(225,30){\circle*{2}}

\put(210,12){\line(0,1){16}} \put(257,12){\line(0,1){16}}
\put(194,12){\line(1,1){16}}

\put(210,7){\makebox(0,0)[t]{\tiny$G\!A$}}
\put(257.5,7){\makebox(0,0)[t]{\tiny$0$}}
\put(194,7){\makebox(0,0)[t]{\tiny$G\!A\!\pl\! 1$}}
\put(210,33){\makebox(0,0)[b]{\tiny$G\!A$}}
\put(257.5,33){\makebox(0,0)[b]{\tiny$0$}}
\put(175,20){\makebox(0,0)[r]{$G\delta_A$}}

\end{picture}
\end{center}

We can prove the following, by proceeding as for the proof of
\LLS-Coherence in Section~2.

\prop{\LLSco-Coherence}{The functor $G$ from \LLSco\ to
$\Delta^{op}$ is faithful.}

We can prove for right strong comonads coherence results analogous
to the two versions of \LRS-Coherence in Section~2.

\section{Coherence for monoidal comonads}
A \emph{monoidal} comonad in a monoidal category $\cal A$ is a
comonad $\langle L,\varepsilon,\delta\rangle$ (in the notation of
\cite{ML98}, Section VI.1) in $\cal A$ together with a natural
transformation whose components are the arrows
\[
\psi_{A,B}\!:LA\otimes LB\str L(A\otimes B),
\]
and together with the arrow $\psi_0\!:I\str LI$, such that $L$
with $\psi$ and $\psi_0$ is a monoidal functor (which means that
we have the equations ${(\psi a)}$, ${(\psi l)}$ and ${(\psi r)}$
with $T$ and $\eta_I$ replaced respectively by $L$ and $\psi_0$),
and we have moreover the equations
\begin{tabbing}
\hspace{1.5em}\=${(\psi
\varepsilon)}$\hspace{3em}\=$\varepsilon_{A\otimes
B}\cirk\psi_{A,B}\,$\=$=\varepsilon_A\otimes\varepsilon_B$,
\\[2ex]
\>${(\psi\delta)}$\>$\delta_{A\otimes B}\cirk\psi_{A,B}$\>$=
L\psi_{A,B}\cirk\psi_{LA,LB}\cirk(\delta_A\otimes\delta_B)$,
\\[3ex]
\>${(\psi_0\epsilon)}$\>$\varepsilon_I\cirk\psi_0\,$\=$=\mj_I$,\\
[2ex]
\>${(\psi_0\delta)}$\>$\delta_I\cirk\psi_0$\>$=L\psi_0\cirk\psi_0$,
\end{tabbing}
which say that $\varepsilon$ and $\delta$ are monoidal natural
transformations. (References concerning the notion of monoidal
comonad are given in Section~1.)

Let \MSco\ be the category of the monoidal comonad freely
generated by an arbitrary set of objects. We define the functor
$G$ from \MSco\ to \Rel\ with the clause for $G\psi_{A,B}$
corresponding to the following picture:
\begin{center}
\begin{picture}(320,50)

\put(142,16){\line(0,1){20}} \put(146,16){\line(0,1){20}}

\put(141.75,16){\line(1,0){4.3}} \put(141.75,36){\line(1,0){4.3}}

\put(170,16){\line(0,1){20}} \put(174,16){\line(0,1){20}}

\put(169.75,16){\line(1,0){4.3}} \put(169.75,36){\line(1,0){4.3}}

\put(132,16){\line(3,2){30}}

\put(132,16){\line(0,1){20}}

\put(130,42){\makebox(0,0)[l]{$L\;A\,\otimes \,L\,B$}}
\put(128.25,8){\makebox(0,0)[l]{$L\,(A\;\otimes\;\,B)$}}

\end{picture}
\end{center}
We also have the clause that says that $G\psi_0$ is the empty
relation from \pr\ to $\{\pr\}$, the clauses for $G\varepsilon_A$
and $G\delta_A$ given above, and the remaining clauses as in
Section~2.

For the proof of \MSco-Coherence we need the following notion of
normal form. An arrow term $f_2\cirk f_1$ of \MSco\ is in normal
form when every factor of the developed arrow term $f_1\!:A\str C$
is an $\varepsilon$-factor or a $\delta$-factor, and $f_2\!:C\str
B$ is an arrow term of the category $\cal M$, which is the free
monoidal category with a single monoidal endo\-functor (see
\cite{DP09a}, Section~3; the family of monoidal endo\-functors of
$\cal M$ is here taken to be the singleton $\{L\}$). It is easy to
see that the equations of \MSco\ yield that every arrow term is
equal to an arrow term in normal form. To ascertain that
$G(f_2\cirk f_1)$ determines uniquely  $Gf_1$, $Gf_2$ and $C$ we
rely on a general proposition about decomposing an arbitrary
binary relation between finite ordinals into three functions.

To formulate this proposition, let $<_l$ be the lexicographical
order on $n\times m$; i.e.\ for $x_1,x_2\in n$ and $y_1,y_2\in m$
we have
\[
(x_1,y_1)<_l(x_2,y_2)\hspace{.5em}\mbox{iff}\hspace{.5em}(x_1<x_2
\hspace{.5em}\mbox{or}\hspace{.5em}(x_1=x_2\hspace{.5em}\mbox{and}
\hspace{.5em}y_1<y_2)).
\]
We call this the \emph{left} lexicographical order, while the
\emph{right} lexicographical order $<_r$ is defined by
\[
(x_1,y_1)<_r(x_2,y_2)\hspace{.5em}\mbox{iff}\hspace{.5em}(y_1<y_2
\hspace{.5em}\mbox{or}\hspace{.5em}(y_1=y_2\hspace{.5em}\mbox{and}
\hspace{.5em}x_1<x_2)).
\]

Let $\langle\nu,\mu,\beta\rangle$ be a triple of functions
$\nu\!:k\str n$, $\mu\!:k\str m$ and $\beta\!:k\str k$, for
$\beta$ a bijection, such that for $z\in k$ and
\begin{tabbing}
\hspace{11em}\=$l_{\nu,\mu,\beta}(z)$ \=
$=_{df}(\nu(z),\mu(\beta(z)))$,
\\[1.5ex]
\>$r_{\nu,\mu,\beta}(z)$\>$=_{df}(\nu(\beta^{-1}(z)),\mu(z))$,
\end{tabbing}
we have for every $u,v\in k$
\[
(\ast)\hspace{1em}\mbox{if}\hspace{.5em}u<v,\hspace{.5em}
\mbox{then}\hspace{.5em}(l_{\nu,\mu,\beta}(u)<_l
l_{\nu,\mu,\beta}(v)\hspace{.5em} \mbox{and}\hspace{.5em}
r_{\nu,\mu,\beta}(u)<_r r_{\nu,\mu,\beta}(v)).
\]

An alternative condition equivalent to $(\ast)$ is to say that
$l_{\nu,\mu,\beta}(z)$ is the $(z\pl 1)$-th pair in the $<_l$
ordering of the image of $l_{\nu,\mu,\beta}$, and analogously with
$r$. The image of $l_{\nu,\mu,\beta}$ coincides with the image of
$r_{\nu,\mu,\beta}$; it coincides also with the set of ordered
pairs $\mu\cirk\beta\cirk\nu^{-1}$, whose cardinality is $k$.

We call triples of functions such as $\langle\nu,\mu,\beta\rangle$
above \emph{coordinated}. That a triple of functions is
coordinated amounts to saying that $\nu$ and $\mu$ are
order-preserving and that in $k$ there are no analogues of the
short circuits and useless crossings of \cite{DP09a} (Section~7).
We can now formulate our general proposition about decomposition.

\prop{Decomposition Proposition}{For every relation $R\subseteq
n\times m$ there is a unique coordinated triple of functions
$\langle\nu,\mu,\beta\rangle$ such that
$R=\mu\cirk\beta\cirk\nu^{-1}$. The domain of $\nu$, $\mu$ and
$\beta$ is the cardinality of $R$.}

\noindent Here $\cirk$ on the right-hand side is composition of
relations, and $\nu^{-1}$ is the relation converse to the function
$\nu$. We denote the cardinality of the set of ordered pairs $R$
by $|R|$. The Decomposition Proposition is illustrated by the
following example:

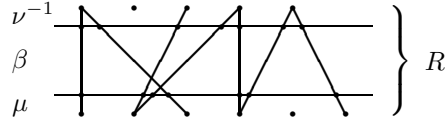
\begin{center}
\begin{picture}(150,40)
{\thicklines \put(20,0){\line(0,1){40}} \put(40,0){\line(1,2){20}}

\put(40,0){\line(1,1){40}} \put(60,0){\line(-1,1){40}}
\put(80,0){\line(0,1){40}} \put(80,0){\line(1,2){20}}
\put(120,0){\line(-1,2){20}}

\put(10,38){\makebox(0,0)[r]{$\nu^{-1}$}}
\put(0,20){\makebox(0,0)[r]{$\beta$}}
\put(0,2){\makebox(0,0)[r]{$\mu$}}}

{\thinlines \put(10,7){\line(1,0){120}}
\put(10,33){\line(1,0){120}}}

\put(20,0){\circle*{2}} \put(40,0){\circle*{2}}
\put(60,0){\circle*{2}} \put(80,0){\circle*{2}}
\put(100,0){\circle*{2}} \put(120,0){\circle*{2}}

\put(20,40){\circle*{2}} \put(40,40){\circle*{2}}
\put(60,40){\circle*{2}} \put(80,40){\circle*{2}}
\put(100,40){\circle*{2}}

\put(20,7){\circle*{2}} \put(43.5,7){\circle*{2}}
\put(47,7){\circle*{2}} \put(53,7){\circle*{2}}
\put(80,7){\circle*{2}} \put(83.5,7){\circle*{2}}
\put(116.5,7){\circle*{2}}

\put(20,33){\circle*{2}} \put(27,33){\circle*{2}}
\put(56.5,33){\circle*{2}} \put(73,33){\circle*{2}}
\put(80,33){\circle*{2}} \put(96.5,33){\circle*{2}}
\put(103.5,33){\circle*{2}}

\put(125,20){\makebox(0,0)[l]{$\displaystyle
\left.\begin{array}{l} \\[3em] \end{array}\right \}$}}

\put(150,20){\makebox(0,0)[l]{$R$}}

\end{picture}
\end{center}
which also makes its truth pretty obvious.

We will however prove this proposition formally. For that, let the
bijection $l_R\!:|R|\str R$ be defined by
\begin{tabbing}
\centerline{$l_R(z)$ is the $(z\pl 1)$-th ordered pair of $R$ in
the ordering $<_l$.}
\end{tabbing}
We define analogously the bijection $r_R\!:|R|\str R$ via $<_r$.

Given $R$, consider the following functions:
\begin{tabbing}
\hspace{11em}\=$\beta_R$ \=$=_{df} r^{-1}_R\cirk l_R$ \=$:|R|\str
|R|$,\kill

\>$\nu_R$ \>$=_{df} \;\;p^1\cirk l_R$ \>$:|R|\str n$,
\\[1.5ex]
\>$\mu_R$ \>$=_{df} \;\,p^2\cirk r_R$ \>$:|R|\str m$,
\\[1.5ex]
\>$\beta_R$ \>$=_{df} r^{-1}_R\cirk l_R$ \=$:|R|\str |R|$,
\end{tabbing}
where $p^1$ and $p^2$ are respectively the first and second
projection with domain $n\times m$.

To show that $l_{\nu_R,\mu_R,\beta_R}$ and
$r_{\nu_R,\mu_R,\beta_R}$ satisfy $(\ast)$ it is enough to verify
that the following lemma holds.

\prop{Lemma~1}{$l_{\nu_R,\mu_R,\beta_R}=l_R$ and
$r_{\nu_R,\mu_R,\beta_R}=r_R$.}

\dkz For the first equation we have
\begin{tabbing}
\hspace{5em}\=$(\nu_R(z),\mu_R(\beta_R(z)))$ \=$=
(p^1(l_R(z)),p^2(r_R(r^{-1}_R(l_R(z)))))$
\\[1.5ex]
\>\>$= (p^1(l_R(z)),p^2(l_R(z)))$
\\[1.5ex]
\>\>$= l_R(z)$,
\end{tabbing}
and analogously for the second equation. \qed

\vspace{2ex}

\noindent Since we also have $R=\mu_R\cirk\beta_R\cirk\nu^{-1}_R$,
the triple $\tau(R)=\langle\nu_R,\mu_R,\beta_R\rangle$ is a
coordinated triple of functions such as required by the
Decomposition Proposition.

To show that $\tau(R)$ is unique we proceed as follows. It is
enough to verify besides $R=\mu_R\cirk\beta_R\cirk\nu^{-1}_R$ that
for every coordinated triple of functions
$\langle\nu,\mu,\beta\rangle$ we have
\[
\tau(\mu\cirk\beta\cirk\nu^{-1})=\langle\nu,\mu,\beta\rangle.
\]
For that we rely on the following lemmata.

\prop{Lemma~2}{$l_{\mu\cirk\beta\cirk\nu^{-1}}=l_{\nu,\mu,\beta}$
and $r_{\mu\cirk\beta\cirk\nu^{-1}}=r_{\nu,\mu,\beta}$.}

\noindent For the proof we rely on the comment after $(\ast)$.

\prop{Lemma~3}{$\nu_{\mu\cirk\beta\cirk\nu^{-1}}=\nu$,
$\mu_{\mu\cirk\beta\cirk\nu^{-1}}=\mu$ and
$\beta_{\mu\cirk\beta\cirk\nu^{-1}}=\beta$.}

\dkz For the first equation we have
\begin{tabbing}
\hspace{5em}\=$\nu_{\mu\cirk\beta\cirk\nu^{-1}}(z)$
\=$=p^1(l_{\mu\cirk\beta\cirk\nu^{-1}}(z))$,
\\[1.5ex]
\>\>$=p^1(\nu(z),\mu(\beta(z)))$,\quad by Lemma~2
\\[1.5ex]
\>\>$=\nu(z)$.
\end{tabbing}
The second equation is derived analogously, while for the third
equation we have
\begin{tabbing}
\hspace{5em}\=$\beta_{\mu\cirk\beta\cirk\nu^{-1}}(z)$
\=$=r^{-1}_{\mu\cirk\beta\cirk\nu^{-1}}(l_{\mu\cirk\beta\cirk\nu^{-1}}(z))$,
\\[1.5ex]
\>\>$=r^{-1}_{\nu,\mu,\beta}(\nu(z),\mu(\beta(z)))$,\quad by
Lemma~2
\\[1.5ex]
\>\>$=\beta(z)$.
\end{tabbing}
since we have
\begin{tabbing}
\hspace{5em}\=$\beta_{\mu\cirk\beta\cirk\nu^{-1}}(z)$ \=\kill
\>$\;\;r_{\nu,\mu,\beta}(\beta(z))$
\>$=\nu(\beta^{-1}(\beta(z)),\mu(\beta(z)))$,
\\[1.5ex]
\>\>$=(\nu(z),\mu(\beta(z)))$.\`$\dashv$
\end{tabbing}

\vspace{1ex}

\noindent This concludes the proof of the Decomposition
Proposition. (The Decomposition Proposition could be used to
obtain a normal form for arrow terms of the category \Rel---a
normal form alternative to the iota normal form of \cite{DP09};
Section 13.)

We can now finish the proof of the following.

\prop{\MSco-Coherence}{The functor $G$ from \MSco\ to \Rel\ is
faithful.}

\dkz We rely on the normal form $f_2\cirk f_1$, which we
introduced before the Decomposition Proposition. In $f_1$ we find
what corresponds to $\nu^{-1}$ in the Decomposition Proposition,
and in $f_2$ what corresponds to $\mu\cirk\beta$. We rely then on
coherence for comonads (see \cite{DP08}, Section~3, and references
therein) and on $\cal M$-Coherence of \cite{DP09a} (Section~4).
\qed

\vspace{2ex}

The normal form of this proof could be refined to $f_4\cirk
f_3\cirk f_2\cirk f_1$ where in the developed arrow term $f_1$
every factor is an $\varepsilon$-factor, in the developed arrow
term $f_2$ every factor is a $\delta$-factor, in the developed
arrow term $f_4$ every factor is a $\psi_0$-factor, and $f_3$ is
an arrow term of $\cal M$ without occurrences of $\psi_0$.

\section{Coherence for symmetric monoidal comonads}
A \emph{symmetric} monoidal comonad is a monoidal comonad in a
symmetric monoidal category whose endo\-functor $L$ is a monoidal
functor (see the beginning of the preceding section) that
satisfies the equation ${(\psi c)}$ of Section~4 with $T$ replaced
by $L$; namely, this endo\-functor is a linear endo\-functor in
the sense of \cite{DP09a} (Section~6).

Let \McSco\ be the category of the symmetric monoidal comonad
freely generated by an arbitrary set of objects. We define the
functor $G$ from \McSco\ to \Rel\ as $G$ from \MSco\ to \Rel, save
that now we have that $GA$ is the number of occurrences of
generating objects, i.e.\ propositional letters, in $A$ plus the
number of occurrences of $L$ in $A$ (see Section~4). We have
moreover a clause for $Gc_{A,B}$ as in Section~4. We can prove the
following.

\prop{\McSco-Coherence}{The functor $G$ from \McSco\ to \Rel\ is
faithful.}

\noindent For the proof we proceed as for \MSco-Coherence, by
relying on a normal form $f_2\cirk f_1$ where $f_1$ is as before
while $f_2$ is an arrow term of the category \Mc, which is the
free symmetric monoidal category with a single linear
endo\-functor; we appeal then to \Mc-Coherence of \cite{DP09a}
(Section~6).

\section{Coherence for cartesian monoidal comonads}
A \emph{cartesian} monoidal comonad is a symmetric monoidal
comonad in a cartesian category, which satisfies moreover the
equation ${(\psi\Delta)}$ of Section~5 with $T$ replaced by $L$.

Let \CSco\ be the category of the cartesian monoidal comonad
freely generated by an arbitrary set of objects. This category may
be taken as axiomatizing identity of deductions in the
$\{\Box,\wedge,\top\}$ fragment of the modal logic S4 (cf.\
\cite{DP08}). We define the functor $G$ from \CSco\ to \Rel\ as
$G$ from \McSco\ to \Rel\ with additional clauses for $G\Delta_A$
and $G\esp_A$ as in Section~5. We can prove the following.

\prop{\CSco-Coherence}{The functor $G$ from \CSco\ to \Rel\ is
faithful.}

\noindent For the proof we proceed as for {$\cal C$}-Coherence in
\cite{DP09a} (Section~8). We rely again on the possibility to
assume that the targets are $\otimes$-free.

\section{Coherence for cocartesian monoidal comonads}
A \emph{cocartesian} monoidal comonad is a symmetric monoidal
comonad in a cocartesian category; we have moreover the equation
${(\psi\;\mbox{\it def}\,)}$ of Section~6. The equation
\begin{tabbing}
\hspace{1.5em}\=${(\psi c)}$\hspace{3em}\=$Tc_{A,B}\cirk
\psi_{A,B}=\psi_{B,A}\cirk c_{TA,TB}$.\kill

\>${(\psi_0\;\mbox{\it def}\,)}$\>$\psi_0=\;!_{LI}$
\end{tabbing}
follows from the assumption that $I$ is an initial object, which
comes with the assumption that we are in a cocartesian category.

More simply, we can define a cocartesian monoidal comonad as a
cocartesian category with a comonad in it. The definitions of
$\psi$ and $\psi_0$ are then given by ${(\psi\;\mbox{\it def}\,)}$
and ${(\psi_0\;\mbox{\it def}\,)}$.

Let \DSco\ be the category of the cocartesian monoidal comonad
freely generated by an arbitrary set of objects. We define the
functor $G$ from \DSco\ to \Rel\ as $G$ from \McSco\ to \Rel\ with
additional clauses for $G\nabla_A$ and $G!_A$ as in Section~6. We
can prove the following.

\prop{\DSco-Coherence}{The functor $G$ from \DSco\ to \Rel\ is
faithful.}

\dkz We rely on a normal form $f_2\cirk f_1$ where every factor of
the developed arrow term $f_1\!:A\str C$ is an
$\varepsilon$-factor or a $\delta$-factor, and $f_2\!:C\str B$ is
an arrow term of the category $\cal D$, which is the free
cocartesian category with a single endo\-functor (see
\cite{DP09a}, Section~9). To ascertain that $G(f_2\cirk f_1)$
determines uniquely $Gf_1$, $Gf_2$ and $C$ we use the
Decomposition Proposition of Section~8. We rely then on coherence
for comonads, as in the proof of \MSco-Coherence in Section~8, and
on $\cal D$-Coherence of \cite{DP09a} (Section~9).\qed

\vspace{4ex}

\noindent {\small {\it Acknowledgement.} Work on this paper was
supported by the Ministry of Science of Serbia (Grants 144013 and
144029). We thank an anonymous referee for finding some typos and
making some useful suggestions to improve the exposition.}

\end{document}